\documentclass[a4paper,11pt,intlimits,oneside]{amsart}
\usepackage{amsfonts, amssymb, amsmath}
\usepackage{enumerate}
\usepackage{esint}
\usepackage{latexsym,amssymb}

\newtheorem{thm}{Theorem}[section]

\newtheorem{cor}[thm]{Corollary}
\newtheorem{lem}[thm]{Lemma}
\newtheorem{prop}[thm]{Proposition}
\newtheorem{defn}[thm]{Definition}
\newtheorem{rem}[thm]{Remark}

\newcommand{\f}{\frac}

\newcommand{\vc}{\infty}
\newcommand{\dx}{d\mu(x)}
\newcommand{\dy}{d\mu(y)}

\newcommand{\dtt}{\f{dt}{t}}
\newcommand{\dyt}{\f{\dy}{V(x,t)}\f{dt}{t}}

\newcommand{\LA}{L_kA^{-1/2}}

\newcommand{\VA}{V^{1/2}A^{-1/2}}
\newcommand{\RR}{X}
\newcommand{\su}{\subset}

\newcommand{\HL}{H^p_{L,w}(X)}

\newcommand{\XL}{X\times (0,\vc)}

\newcounter{rek}
\title[Weighted Hardy spaces associated to operators]{ Weighted Hardy spaces associated to operators and boundedness of singular integrals}         

\author{The Anh Bui \& Xuan Thinh Duong}


\begin{document}
\maketitle

\begin{abstract}
Let $(X, d, \mu)$ be a space of homogeneous type, i.e. the measure $\mu$ satisfies doubling (volume)
property with respect to the balls defined by the metric $d$.
Let $L$ be a non-negative self-adjoint
operator on $L^2(X)$. Assume that the semigroup of $L$ satisfies the Davies-Gaffney estimates.
In this paper, we study the weighted Hardy spaces $H^p_{L,w}(X)$, $0 < p \le 1$, associated to the
operator $L$ on the space $X$. We establish the atomic and the molecular
characterizations of elements in $H^p_{L,w}(X)$. As applications, we obtain the boundedness on $\HL$ for
 the generalized Riesz transforms associated to $L$ and for the spectral multipliers of $L$.
\end{abstract}

 \tableofcontents

\section{Introduction}
The theory of Hardy spaces has been a central part of modern harmonic analysis.
While the classical Hardy spaces on $\mathbb R^n$ can be characterized
by certain estimates via the Laplacian, the study  on the Hardy spaces associated to operators has been
intensive recntly, see for example \cite{ADM, DY2, HM, HLMMY} and their
references. Here we will give only a brief account of some  recent studies.
 In \cite{ADM},  the Hardy spaces $H^1_{L}$
associated to an operator $L$ was introduced and studied under the
assumption that the  heat kernel of $L$ satisfies a pointwise
Poisson upper bound. The BMO space associated to such an $L$
was introduced in \cite{DY1}  and it was shown  in \cite{DY2} that the BMO space associated to $L$
is the dual space of the Hardy space $H^1_{L^*}$ associated to the
adjoint operator $L^*$. Recently  the Hardy space $H^1$ associated to Hodge Laplacian on
a Riemannian manifold was studied in \cite{AMR}. Meanwhile the Hardy spaces associated to a second order
divergence form elliptic operator $L$ on $\mathbb{R}^n$ with complex
coefficients was investigated in \cite{HM}.
The study of the Hardy spaces $H^p_L(X) , 1\leq p<\infty,$ on a metric
space $X$ associated to a non-negative self adjoint operator $L$
satisfying Davies-Gaffney estimates was carried out in   \cite{HLMMY}.\\

It is natural to study   weighted Hardy
space $H^p_{L,w} , 1\leq p<\infty$ associated to an operator $L$ with an appropriate weight $w$.
It is known  that the  classical weighted Hardy
spaces $H^p_w(\mathbb R^n)$ can be considered as weighted Hardy spaces
associated to the Laplacian. See, for example,  \cite{G}, \cite{ST}
and their references. \\

This paper is inspired by the recent work of Song and Yan \cite{SY}
in which they introduced the weighted Hardy spaces
$H^1_{L,w}(\mathbb R^n)$ associated to an operator $L$ initially
defined on $L^2(\mathbb R^n)$ with the assumptions that $L$ is
non-negative self-adjoint on $L^2$ (assumption (H1) in Section 3.1)
and that $L$ has Gaussian heat kernel bounds (assumption (H3) in
Section 3.1). As an application, it was shown in \cite{SY} that when
the operator $L$ is the Schr\"odinger operator with a non-negative
potential, the Riesz transform associated to $L$ is bounded from the
weighted Hardy space $H^1_{L,w}(\mathbb R^n)$  to the classical
weighted Hardy space $H^1_w(\mathbb R^n)$.\\

In this paper, we extend the work in \cite{SY} on $H^1_{L,w}(\mathbb R^n)$ in several aspects.
We define and study weighted Hardy space $H^p_{L,w}(X)$  in the general setting:

\medskip

(i) The underlying space $X$ is a space of homogeneous type, i.e. a metric space  with doubling property.

\medskip

(ii) The index $p$ is in the range  $0 < p \le1$.

\medskip

(iii) $L$ is assumed to satisfy the weaker assumption of  Davies-Gaffney estimate (assumption (H2) in Section 3.1) instead of
 the stronger assumption of Gaussian heat kernel bound.\\

In comparison with \cite{SY}, our approach is different from that in \cite{SY}. Let us remind that  the Gaussian
heat kernel upper bound and the setting of Euclidean space $\mathbb{R}^n$ seem to be indispensable and play an essential role in the approach of \cite{SY}.
One of the key estimates used in \cite{SY} is the equivalences of different weighted
area integral norms in \cite{ST}. Moreover, their approach relied heavily on a geometric argument. However, in the setting of homogeneous spaces without Gaussian upper bound condition, it is not clear whether or not their approach still work. In this paper, using the different approach by introducing new
weighted tent spaces (see Section 3.3) and exploiting the similar approach to that of \cite{HLMMY}, we extend the results in \cite{SY} to spaces of homogeneous type $X$ under the weaker assumption of  Davies-Gaffney estimate.\\

To demonstrate practical applications of  our study of weighted Hardy
spaces, we consider certain singular integral operators whose
kernels are not smooth enough for the operators to belong to the class of
standard Calder\'on-Zygmund operators. We show that
  some of these singular integrals are bounded on weighted Hardy spaces
 associated to an (appropriate) operator $L$.\\

 The layout of this paper is as follows. In Section 2, we review the concept of doubling space
 and the main  properties of the   Muckenhoupt weights and  reverse H\"older weights.
 In Section 3, we  introduce the weighted Hardy spaces associated to operators $H^p_{L,w}(X)$ for $0 < p \le 1$
 by using the area integral norm, then
obtain an atomic characterization of elements in $H^p_{L,w}(X)$.
 In Section 4, we show that  certain singular integral operators are bounded on $H^p_{L,w}(X)$.
 These operators have non-smooth kernels and they include the Riesz transforms associated to magnetic Schr\"odinger
operators,  Riesz transforms associated to  an operator and spectral multipliers of  non-negative self-adjoint operator.\\

\medskip

After finishing our paper, we had learned that in \cite{YY} the authors studied Musielak-Orlicz Hardy spaces associated to such an operator $L$. However, it seems that there are some differences in obtaining the atomic decomposition of the weighted tent spaces and the condition on the weights between our paper and \cite{YY}. Moreover, the the kind of singular integrals considered in our paper is also different from that in \cite{YY}. So, the results obtained in this paper are still interesting in their own rights.

\section{Preliminaries}

\subsection{Doubling metric spaces}

Let $X$ be a metric space, with distance $d$ and $\mu$ is a
nonnegative, Borel, doubling measure on $X$. Throughout this paper,
we assume that $\mu (X)=\infty$.\\
Denote by $B(x, r)$ the open ball of radius $r >0$ and center $x\in
M$, and by $V (x, r)$ its measure $\mu(B(x, r))$. The doubling
property of $\mu$ provides that there exists a constant $C>0$ so
that
\begin{equation}\label{doublingpro}
V(x,2r)\leq CV(x,r)
\end{equation}
for all $x\in X$ and $r>0$.\\
Notice that the doubling property (\ref{doublingpro}) implies that
following property that
\begin{equation}\label{doublingpro1}
V(x,\lambda r)\leq C\lambda^nV(x,r),
\end{equation}
for some positive constant $n$ uniformly for all $\lambda\geq 1,
x\in M$ and $r>0$. There also exists a constant $0\leq N\leq n$ such
that
\begin{equation}\label{doublingpro2}
V(x,r)\leq C\Big(1+\frac{d(x,y)}{r}\Big)^NV(y,r),
\end{equation}
uniformly for all $x,y\in X$ and $r>0$.\\

To simplify notation, we will often just use $B$ for $B(x_B, r_B)$ and $V(E)$ for $\mu(E)$ for any measurable subset $E\subset X$.
Also given $\lambda > 0$, we will write $\lambda B$ for the
$\lambda$-dilated ball, which is the ball with the same center as
$B$ and with radius $r_{\lambda B} = \lambda r_B$. For each ball
$B\subset X$ we set
$$
S_0(B)=B \ \text{and} \ S_j(B) = 2^jB\backslash 2^{j-1}B \
\text{for} \ j\in \mathbb{N}.
$$

\subsection{Muckenhoupt weights}

Throughout this article, we shall denote $w(E) :=\int_E w(x)\dx$ for
any measurable set $E \subset X$. For $1 \leq p \leq \infty$ let $p'$
be the conjugate exponent of $p$, i.e. $1/p + 1/p' = 1$.

We first introduce some notation. We use the notation
$$
\fint_B h(x)\dx=\f{1}{V(B)}\int_Bh(x)\dx.
$$
A weight $w$ is a non-negative measurable and locally integrable function on $X$.
We say that $w \in A_p$, $1 < p < \infty$, if there exists a
constant $C$ such that for every ball $B \subset X$,
$$
\Big(\fint_B w(x)\dx\Big)\Big(\fint_B w^{-1/(p-1)}(x)\dx\Big)^{p-1}\leq C.
$$
For $p = 1$, we say that $w \in A_1$ if there is a constant $C$ such
that for every ball $B \subset X$,
$$
\fint_B w(x)\dx \leq Cw(x) \ \text{for a.e. $x\in B$}.
$$
We set $A_\vc=\cup_{p\geq 1}A_p$.

The reverse H\"older classes are defined in the following way: $w
\in RH_q, 1 < q < \infty$, if there is a constant $C$ such that for
any ball $B \subset X$,
$$
\Big(\fint_B w^q(y) \dy\Big)^{1/q} \leq C \fint_B w\dx.
$$
The endpoint $q = \infty$ is given by the condition: $w \in
RH_\infty$ whenever, there is a constant $C$ such that for any ball
$B \subset X$,
$$
w(x)\leq C \fint_B w(y)\dy  \ \text{for a.e. $x\in B$}.
$$
Let $w \in A_\vc$, for $1\leq p <\infty$, the weighted spaces $L^p_w(X)$
can be defined by
$$\Big\{f :\int_{X} |f(x)|^p w(x)\dx < \infty\Big\}$$
with the norm
$$\|f\|_{L^p_w(X)}=\Big(\int_{X} |f(x)|^p w(x)\dx\Big)^{1/p}$$.\\

We sum up some of the properties of $A_p$ classes in the following results, see \cite{Du}.
\begin{lem}\label{weightedlemma1}
The following properties hold:
\begin{enumerate}[(i)]
\item $A_1\subset A_p\subset A_q$ for $1\leq p\leq q\leq \infty$.
\item $RH_\infty \su RH_q \su RH_p$ for $1\leq p\leq q\leq \infty$.
\item If $w \in A_p, 1 < p < \vc$, then there exists $1 < q < p$ such that $w \in A_q$.
\item If $w \in RH_q, 1 < q < \vc$, then there exists $q < p < 1$ such that $w \in RH_p$.
\item $A_\vc =\cup_{1\leq p<\vc}A_p \subset \cup_{1< p\leq \vc}RH_p$
\end{enumerate}
\end{lem}
\begin{lem}\label{weightedlemma2}
For any ball $B$, any measurable subset $E$ of $B$ and $w \in A_p, p \geq  1$,
 there exists a constant $C_1 > 0$ such that
$$
C_1\Big(\f{V(E)}{V(B)}\Big)^p\leq \f{w(E)}{w(B)}.
$$
If $w \in RH_r, r > 1$. Then, there exists a constant $C_2 > 0$ such that
$$
\f{w(E)}{w(B)}\leq C_2\Big(\f{V(E)}{V(B)}\Big)^{\f{r-1}{r}}.
$$
\end{lem}
From the first inequality of Lemma \ref{weightedlemma2}, if $w\in
A_1$ then there exists a constant $C>0$ so that for any ball $B\subset X$ and $\lambda >1$, we have
$$
w(\lambda B)\leq C\lambda^n w(B).
$$

\section{Weighted Hardy spaces associated to operators}

\subsection{Definition of weighted Hardy spaces}
In this paper we consider the following conditions:

\smallskip

\noindent
 ${\bf(H 1)}$ $L$ is a non-negative self-adjoint operator on $L^2({X})$;

\smallskip

\noindent ${\bf(H 2)}$ The operator $L$ generates an analytic
semigroup $\{e^{-tL}\}_{t>0}$ which satisfies the Davies-Gaffney
condition. That is, there exist constants $C$, $c>0$ such that for
any open subsets $U_1,\,U_2\subset X$,

\begin{equation}\label{D-Gestimate}
|\langle e^{-tL}f_1, f_2\rangle| \leq C\exp\Big(-{{\rm
dist}(U_1,U_2)^2\over c\,t}\Big)
\|f_1\|_{L^2(X)}\|f_2\|_{L^2(X)},\quad\forall\,t>0,
\end{equation}

\noindent for every $f_i\in L^2(X)$ with $\mbox{supp}\,f_i\subset
U_i$, $i=1,2$, where ${\rm dist}(U_1,U_2):=\inf_{x\in U_1, y\in U_2}
d(x,y)$.

 \medskip

${\bf(H3)}$ The kernel of $e^{-tL}$ denote by $p_t(x,y) $ which
satisfies the Gaussian upper bound. That is, there exist constants
$C$, $c>0$ such that for almost every $x,y \in X$,
\begin{equation}\label{Gaussianub}
|p_t(x,y)|\leq \f{C}{V(x,\sqrt{t})}\exp\Big(-\f{d(x,y)^2}{ct}\Big),
\forall t>0.
\end{equation}

 It is not difficult to show that condition $(H3)$ implies $(H2)$.\\

Let $L$ be an operator satisfying $(H1)$ and $(H2)$. Set
$$
H^2(X):=\overline{\mathcal{R}(L)}=\{Lu\in L^2(X): u\in \mathcal{D}(L)\}
$$
where $\mathcal{D}(L)$ is the domain of $L$.

It is known that $L^2(X)=\overline{\mathcal{R}(L)}\oplus
\mathcal{N}(L)$, where $\mathcal{R}(L)$ and $\mathcal{N}(L)$ stand
for the range and the kernel of $L$, and the sum is orthogonal.

For $w\in A_\vc$ and $0<p<\vc$, the Hardy space $H_{L,w}^p(X)$ is
defined as the completion of
$$
\{f\in \overline{\mathcal{R}(L)}:\|S_Lf\|\in L^p(w)\}
$$
in the norm $\|f\|_{\HL}=\|S_Lf\|_{L^p(w)}$, where
$$
S_Lf(x)=\Big(\int\int_{d(x,y)<t}|t^2Le^{-t^2L}f(y)|^2\dyt\Big).
$$

\begin{rem}
When $w=1$ the Hardy spaces $\HL$ were introduced in \cite{HLMMY}
and $p=1$. The particular case when $p=1$, $X=\mathbb{R}^n$, and $L$
satisfies (H1) and (H3) the Hardy space $H^1_{L,w}(X)$ was studied
in \cite{SY}.
\end{rem}

We next describe the notion of an $(M, p, w)$-atom and an $(M, p,
w,\epsilon)$-molecule associated to the operator $L$.

\begin{defn}
Suppose $w \in A_\vc$. We say that a function $a\in L^2(X)$ is
an $(M, p, w)$-atom associated to an operator $L$, if there
exists a function $b$ which belongs to $\mathcal{D}(L^M)$, the domain of $L^M$, and a ball
$B$ of $X$ such that
\begin{enumerate}[(i)]
\item $a=L^Mb$;
\item {\rm supp}\,$L^kb\su B, k=0,1,\ldots, M$;
\item $\|(r_B^2L)^kb\|_{L^2(X)}\leq r_B^{2M}V(B)^{1/2}w(B)^{-1/p}, k=0,1,\ldots, M$.
\end{enumerate}
\end{defn}

\begin{defn}
Suppose $w \in A_\vc$. We say that a function $m\in L^2(X)$ is
an $(M, p, w,\epsilon)$-molecule associated to an operator $L$, if
there exists a function $b$ which belongs to $\mathcal{D}(L^M)$, and a
ball $B$ of $X$ such that
\begin{enumerate}[(i)]
\item $a=L^Mb$;
\item $\|(r_B^2L)^kb\|_{L^2(S_j(B))}\leq 2^{-j\epsilon}r_B^{2M}V(2^jB)^{1/2}w(2^jB)^{-1/p}, k=0,1,\ldots,
M$ and $j=0,1,\ldots$.
\end{enumerate}
\end{defn}

It is not difficult to show that an $(M, p, w)$-atom associated to
the ball $B$ is also an $(M, p, w,\epsilon)$-atom associated to the
ball same ball $B$.\\

We will say that $f =\sum_j \lambda_ja_j$ is an atomic $(M, p,
w)$-representation if $\{\lambda_j\}_{j=0}^\vc\in l^p$, each $a_j$ is a
$(M, p, w)$-atom, and the sum converges in $L^2(X)$. Set
$$
\mathbb{H}^p_{L, M, w, at}(X) := \{f : f \ \text{has an atomic $(M, p,
w)$-representation}\}
$$
with the norm
$$\|f\|_{\mathbb{H}^p_{L, M, w, at}(X)}=\inf\{\Big(\sum_{j=0}^\vc|\lambda_j|^p\Big)^{1/p}: f =\sum_j \lambda_ja_j \ \text{is an atomic $(M, p,
w)$-representation}\}.
$$
and we also define the Hardy space $H^p_{L, M, w, at}(X)$ as the completion of $\mathbb{H}^p_{L, M, w, at}(X)$ in norm $\|\cdot\|_{\mathbb{H}^1_{L, M, w, at}(X)}$.

\subsection{Finite  propagation speed for the wave equation}
Let $L$ satisfy ($H1$) and ($H2$) and let $K_{\cos(t\sqrt{L})}$ be
the kernel of the operator $\cos(t\sqrt{L})$. Then there exists a
constant $c_0 > 0$ such that
\begin{equation}\label{finitepropagation}
{\rm supp}\,K_{\cos(t\sqrt{L})}\subset \{(x,y)\in X\times X:
d(x,y)\leq c_0t\},
\end{equation}
see for example \cite{S}.

\begin{lem}\label{lem1 finitepropagation}
Let $\varphi\in C^\vc_0(\mathbb{R})$ be even and {\rm supp}\,$\varphi\subset (-c_0^{-1}, c_0^{-1})$, where $c_0$ as in (\ref{finitepropagation}). Let $\Phi$ denote the Fourier transform of $\varphi$. Then for every $k\in \mathbb{N}$ and $t>0$, the kernel $K_{(t^2L)^k\Phi(t\sqrt{L})}$ of $(t^2L)^k\Phi(t\sqrt{L})$ satisfies
$$
{\rm supp}\,K_{(t^2L)^k\Phi(t\sqrt{L})}\subset \{(x,y)\in X\times X:
d(x,y)\leq t\}.
$$
\end{lem}

\begin{lem}
Let $L$ be an operator satisfying (H1) and (H2). For every $k = 0,
1, \ldots $, the operator $(tL)^ke^{-tL}$ satisfies Davies-Gaffney
estimates (H2).
\end{lem}

For the proof we refer the reader to Lemma 3.5 in \cite{HLMMY}.\\

\subsection{Weighted tent spaces}
For the measurable function $F$ defined on $X\times (0,\vc)$, we set
$$
\mathcal{A}(F)(x)=\Big(\int_{\Gamma(x)}|F(y,t)|^2 \dyt\Big)^{1/2}
$$
where $\Gamma(x)$ is the cone $\{(y,t): d(x,y)<t\}$.

For $0<p\leq 1$ and $w\in L^1_{loc}(X)$ we introduce the tent space
$T^p_w(X)$ as those functions $F$ such that $\mathcal{A}F\in L^p_(X)$
and we set $\|F\|_{T^p_w(X)}=\|\mathcal{A}(F)\|_{L^p_w(X)}$. Note that our weighted tent spaces $T^p_w(X)$ can be
considered an extension of those in \cite{CMS} when $w\equiv 1$ and
$X=\mathbb{R}^n$. In the case that $w\equiv 1$, we write $T^p(X)$
instead of $T^p_w(X)$. Another version of weighted tent spaces was
used in \cite{HSV} but this version is not suitable to our
purpose.

Given the ball $B$ we denote by $\widehat{B}$ the tent over $B$,
i.e, $\widehat{B}=\{(y,t): d(x,y)+t<r\}$. Now a function $a(y,t)$ is
said to be an $T^p_w(X)$-atom whenever it is supported in
$\widehat{B}$ and
\begin{equation}\label{defn-atomictentspace}
\Big(\int_{\widehat{B}}|a(y,t)|^2\f{w(B(y,t))}{V(y,t)}\f{\dy
dt}{t}\Big)^{1/2}\leq w(B)^{\f{1}{2}-\f{1}{p}}.
\end{equation}
Let $W(y,t)=\f{w(B(y,t))}{V(y,t)}$. Then the LHS of (\ref{defn-atomictentspace}) is just
$\|a\|_{L^2(W)}$. It is not difficult to check that an $T^p_w(X)$-atom
belongs to $T^p_w(X)$ for $w\in A_1\cap RH_{\f{2}{2-p}}$.

An important result concerning weighted tent spaces is that
 each function in $T^p_w(X)$ has an atomic decomposition. More
 precisely, we have the following result.
\begin{thm}\label{atomicdecomposition}
Let $w\in A_1\cap RH_{\f{2}{2-p}}$ and $F\in T^p_w(X), 0<p\leq 1$. Then there exist a
sequence of $T^p_w(X)$-atoms $\{a_j\}_{j}$ and the sequence of numbers
$\{\lambda_j\}_{j}$ such that
\begin{equation}\label{convergenceinatomicdecomposition}
F=\sum_{j}\lambda_j a_j
\end{equation}
and
\begin{equation}\label{eq2-atomicdecomposetentspaces}
\sum_{j}|\lambda_j|^p \leq C\|F\|_{T^p_w(X)}^p.
\end{equation}

Moreover, if $F\in T^p_w(X)\cap T^2(X)$ then the series in (\ref{convergenceinatomicdecomposition}) converges in both $T^p_w(X)$ and $T^2(X)$.
\end{thm}

Before giving the proof for Theorem \ref{atomicdecomposition} we
need the following technical lemma which is an extension of
\cite[Proposition 3]{HSV} to spaces of homogeneous type.
\begin{lem}\label{lem1-atomicdecom}
Let $w\in A_\vc$ and $B$ a ball in $X$. Then there exists a constant
$C$ such that, for every measure function $F$ defined on $\XL$ and
every measurable set $E\subset B$, we have
$$
\int_{\widehat{B}\backslash
\widehat{\Omega}}|F(y,t)|^2\f{w(B(y,t))}{V(y,t)}\f{\dy dt}{t}\leq C
\int_{B\backslash E}|\mathcal{A}(F)(x)|^2w(x)\dx,
$$
where $\Omega=\{x\in B: \mathcal{M}(\chi_E)(x)>\gamma\}$ for $\gamma\in (0,1)$,
$\gamma$ sufficiently small ($\mathcal{M}$ is the Hardy-Littlewood maximal function).
\end{lem}

\emph{Proof:} We adapt the ideas in \cite[Proposition 3]{HSV} to our situation. Set $S=\{(x,y,t)\in (B\backslash E)\times
(\widehat{B}\backslash\widehat{\Omega}): y\in \Gamma(x)\}$ and
$S(y,t)=B(y,t)\cap (B\backslash E)$. For $(y,t)\in
\widehat{B}\backslash\widehat{\Omega}$ we can pick $x_0$ such that
$x_0\notin \Omega$ and $x_0\in B(y,t)\subset B$. This implies
$\mathcal{M}(\chi_E)(x_0)\leq \gamma$, and hence
$$
\mu(B(y,t)\cap E)\leq \gamma V(y,t).
$$
Since $w\in A_\vc$, $w\in RH_r$ for some $1<r<\vc$. This together with Lemma \ref{weightedlemma2} implies $w(B(y,t)\cap E)\leq
C_\gamma w(B(y,t))$.

This together with $B(y,t)\subset B$ gives
$$
w(B(y,t)\cap (B\backslash E))> (1-C_\gamma) w(B(y,t))
$$
provided with sufficiently small $\gamma$.\\

Then we have
\begin{equation*}
\begin{aligned}
\int_{\widehat{B}\backslash
\widehat{\Omega}}|F(y,t)|^2\f{w(B(y,t))}{V(y,t)}\f{\dy dt}{t}&\leq
C\int_{S}|F(y,t)|^2\int_{S(y,t)}w(x)\dx \f{\dy
}{V(y,t)}\f{dt}{t}\\
&= C\int_{\widehat{B}\backslash
\widehat{\Omega}}|F(y,t)|^2w(x)\dx \f{\dy }{V(y,t)}\f{dt}{t}\\
&\leq C\int_{B\backslash
E}w(x)\Big(\int_{\Gamma(x)}|F(y,t)|^2\f{\dy }{V(y,t)}\f{dt}{t}\Big)\dx\\
&= C\int_{B\backslash E}|\mathcal{A}(F)(x)|^2w(x)\dx.
\end{aligned}
\end{equation*}
This completes our proof.\\

We now prove Theorem \ref{atomicdecomposition}.\\

\emph{Proof of Theorem \ref{atomicdecomposition}:} To prove this theorem, we exploit the standard arguments, see for example \cite{CMS, HSV, R}.

For $k\in
\mathbb{Z}$ let $E_k=\{x: \mathcal{A}F(x)>2^k\}$ and $\Omega_k=\{x:
\mathcal{M}(\chi_{E_k})(x)>\gamma\}$ for some $\gamma\in (0,1)$. Then
$E_k\subset \Omega_k$ and $\mu(\Omega_k)\leq C\mu(E_k)$ for all $k$.
Moreover since $w\in A_1$ we also obtain $w(\Omega_k)\leq
C(w)w(E_k)$ for sufficiently small $\gamma$. It can be verified that
supp\,$F\subset \cup \widehat{\Omega}_k$.

For each $k$, due to \cite[Theorem 1.3, Chapter III]{CW} we can pick
a family of balls $\{Q_k^j\}_{j}$ of $\Omega_k$ satisfying the
following three conditions:
\begin{enumerate}[(i)]
\item $\Omega_k = \cup_j Q_k^j$;
\item there exists a constant $\kappa$ which depends only of $X$ such that $\sum_j \chi_{Q_k^j}\leq
\kappa$;
\item there exists a constant $C_0$ such that $C_0Q_k^j\cap (\Omega_k)^c\neq
\emptyset$.
\end{enumerate}

Taking $C_1>C_0+1$ and setting $B_k^j=C_1Q_k^j$ we have
$$
\widehat{\Omega}_k\backslash \widehat{\Omega}_{k+1}\subset \cup_j
A^j_k
$$
with
$$
A_k^j=\widehat{B}_k^j\cap (Q_k^j\times (0,\vc))\cap
(\widehat{\Omega}_k\backslash \widehat{\Omega}_{k+1}).
$$
We define $a_k^j=2^{-(k+1)}w(B_k^j)^{-1/p}F\chi_{A_k^j}$ and
$\lambda_k^j= 2^{(k+1)}w(B_k^j)^{1/p}$. Then obviously,
$$
F=\sum_{k,j}\lambda_{k}^j a_{k}^j.
$$
Moreover, we have,
 by Lemma \ref{lem1-atomicdecom},
\begin{equation*}
\begin{aligned}
\|a_k^j\|^2_{L^2(W)}&\leq
2^{-2(k+1)}w(B_k^j)^{-2/p}\int_{\widehat{B}_k^j\backslash
\widehat{\Omega}_{k+1}}|F(y,t)|^2\f{w(B(y,t))}{V(y,t)}\f{dydt}{t}\\
&\leq C 2^{-2(k+1)}w(B_k^j)^{-2/p}\int_{B_k^j\backslash E_{k+1}}|\mathcal{A}(F)(x)|^2w(x)\dx\\
&\leq C w(B_k^j)^{1-2/p}.
\end{aligned}
\end{equation*}
Therefore, $a_k^j$ is a multiple $T^p_w(X)$ atom by a
constant.

Furthermore, we have
\begin{equation*}
\begin{aligned}
\sum_{k,j}|\lambda_{k}^j|^p&=\sum_{k,j}2^{p(k+1)}w(B_k^j)\\
&\leq C\sum_{k,j}2^{p(k+1)}w(Q_k^j)\\
&\leq C\sum_{k,j}2^{p(k+1)}w(Q_k^j).
\end{aligned}
\end{equation*}
Since $\{Q_k^j\}_j$ is finite overlap, we have
\begin{equation*}
\begin{aligned}
\sum_{k,j}|\lambda_{k}^j|^p
&\leq C\sum_{k}2^{p(k+1)}w(\Omega_k)\\
&\leq C\sum_{k}2^{p(k+1)}w(E_k)\\
& =C \sum_{k}2^{p(k+1)}w\{x: \mathcal{A}F(x)>2^k\}\\
&\leq C\|\mathcal{A}F\|_{L^p_w(X)}^p=\|F\|_{T^p_w(X)}^p.
\end{aligned}
\end{equation*}

At this stage, the same argument as in \cite[Proposition 3.1]{JY} shows that if
$F\in T^p_w(X)\cap T^2(X)$ then the identity
(\ref{convergenceinatomicdecomposition}) converges in both $T^p_w(X)$
and $T^2(X)$.

Let us denote by $T^{p,c}_w(X)$ and $T^{p,c}(X)$ the spaces of those
functions in $T^p_w(X)$ and $T^p(X)$ with bounded support
respectively. Here, by a function $f$ on $X \times (0,\vc)$ with bounded
support, we shall mean that there exist a ball $B\subset  X$ and $0 < c_1 < c_2 < \vc$ such that supp\,$f \subset B\times (c_1, c_2)$.

\begin{lem}
For $w\in A_1\cap RH_{\f{2}{2-p}}$, the space $T^{p,c}_w(X)$ and $T^{2,c}(X)$
coincide for all $p\in (0,1]$.
\end{lem}

\emph{Proof:} The same argument as in \cite[p. 306]{CMS} gives
$T^{p,c}(X)\subset T^{2,c}(X)$ for all $p\in (0,2]$. By the H\"older
inequality we obtain $T^{2,c}(X)\subset T^{p,c}(X)$ for all $p\in
(0,2]$, and hence $T^{p,c}(X)$ and $T^{2,c}(X)$ coincide for all $p\in
(0,2]$.

Suppose that $f\in T^{2,c}(X)$ with supp\,$f\subset K$ for some bounded
set $K$. Then there exists a ball $\widehat{B}\supset K$ such that
supp\,$\mathcal{A}f \subset B$. By the H\"older inequality we have for all
$p\in (0,1]$
\begin{equation*}
\begin{aligned}
\|\mathcal{A}f\|_{L^p_w(X)}^p&=\int_B|\mathcal{A}f(x)|^pw(x)\dx\\
&\leq \Big(\int_B|\mathcal{A}f(x)|^2\dx\Big)^{p/2}\Big(\int_B
w^{\f{2}{2-p}}\Big)^{\f{2-p}{p}}.
\end{aligned}
\end{equation*}
Since $w\in A_1\cap RH_{\f{2}{2-p}}$ there exists $c(p,w,B)$ such that
$$
\Big(\int_B w^{\f{2}{2-p}}\Big)^{\f{2-p}{p}}\leq c(p,w,B).
$$
This implies $\|\mathcal{A}f\|_{L^p_w(X)}^p \leq c\|f\|^p_{T^2(X)}$.
Therefore, $T^{2,c}(X)\subset T^{p,c}_w(X)$ for all $p\in (0,1]$.

Conversely, for any $f\in T^{p,c}_w(X)$ with supp\,$f\subset K$ for some
bounded set $K$. Let $B$ be the ball satisfying $K\subset
\widehat{B}$ and supp\,$\mathcal{A}f\subset B$. For any $0<r<p$ we
have
\begin{equation*}
\begin{aligned}
\|f\|^r_{T^r(X)}=\|\mathcal{A}f\|_{L^r(X)}^r&=\int_B
|\mathcal{A}f(x)|^r\dx=\int_B
|\mathcal{A}f(x)|^rw(x)^{r/p}w(x)^{-r/p}\dx\\
&\leq \Big(\int_B
|\mathcal{A}f(x)|^pw(x)\Big)^{r/p}\Big(\int_Bw(x)^{(-r/p)(p/r)'}\dx\Big)^{\f{1}{(p/r)'}}.
\end{aligned}
\end{equation*}
Since $w\in A_1$, $w^{(-r/p)(p/r)'}\in A_\vc$. This implies
$w^{(-r/p)(p/r)'}$ is a doubling measure. So,
$$
\Big(\int_Bw(x)^{(-r/p)(p/r)'}\dx\Big)^{\f{1}{(p/r)'}}\leq C.
$$
It therefore follows $T^{p,c}_w(X)\subset  T^{r,c}(X)$. At this stage,
using the fact that $T^{r,c}(X)$ and $T^{2,c}(X)$ coincide, the proof is
complete.

\subsection{Atomic characterization of weighted Hardy spaces $\HL$}

Let $\varphi$ and $\Phi$ be as in Lemma \ref{lem1
finitepropagation}. Setting, $\Psi(t)=t^{2(M+1)}\Phi(t), t\in
\mathbb{R}_+$, then for all $f \in L^{2}_c(\XL)$, the set of all
functions in $L^{2}(\XL)$ with bounded support, and for $x \in X$,
define
$$
\pi_{\Psi,L}f=c_{\Psi}\int_0^\vc\Psi(t\sqrt{L})(f(\cdot,t))(x)\dtt,
$$
where $c_{\Psi}$ is a constant such that
$$
1=c_{\Psi}\int_0^\vc\Psi(t)t^2e^{-t^2}\dtt.
$$
\begin{prop}\label{Pro1-mappingfrom tentspace to Hardyspace}
Let $L$ satisfy (H1) and (H2), $M>\f{n(2-p)}{4p}$, $p\in (0,1]$ and $w\in A_1\cap RH_{\f{2}{2-p}}$. Then
\begin{enumerate}[(i)]
\item the operator $\pi_{\Psi,L}$, initially defined on $T^{2,c}(X)$, extends
to a bounded linear operator from $T^2(X)$ to $L^2(X)$.

\item the operator $\pi_{\Psi,L}$, initially defined on $T^{p,c}_w(X)$, extends
to a bounded linear operator from $T^p_w(X)$ to $H^p_{L,w}(X)$.
\end{enumerate}
\end{prop}
\emph{Proof:} (i) The proof of (i) is standard and we omit the
detail.

(ii) Let $f\in T^{p,c}_w(X)$. It is easy to see that $f\in T^{2,c}(X)$.
Then by Theorem \ref{atomicdecomposition} we have
$$
\pi_{\Psi,L}f=\sum_{j=1}^\vc \lambda_j \pi_{\Psi,L}(a_j)
$$
in $L^2(X)$ with $\{\lambda_j\}$ and $\{a_j\}$ satisfying
(\ref{convergenceinatomicdecomposition}) and
(\ref{eq2-atomicdecomposetentspaces}). Since $S_L$ is bounded on
$L^2(X)$, we obtain that
$$
\|S_L(\pi_{\Psi,L}f)\|^p_{L^p_w(X)}\leq
\sum_{j}|\lambda_j|^p\|S_L(\pi_{\Psi,L}a_j)\|^p_{L^p_w(X)}.
$$
We now claim that $\pi_{\Psi,L}a_j$ is a multiple of an
$(M,p,w)$ atom for each $j$. Indeed, we can write
$\pi_{\Psi,L}a_j=L^Mb_j$ where
$$
b_j=c_{\Psi}\int_0^\vc t^{2M}t^2L\Phi(t\sqrt{L})(a_j(\cdot,t))\dtt.
$$
Let us note that for each $j$ there exists some ball $B_j$ such
that supp\,$a_j\subset \widehat{B_j}$. Therefore, by Lemma \ref{lem1
finitepropagation} we have supp\,$L^kb_j\subset  B_j$ for all
$k=0,\ldots, M$. On the other hand, for any $h\in L^2(B_j)$, by the
H\"older inequality we have
\begin{equation*}
\begin{aligned}
\Big|\int &(r^2_{B_j}L)^kb_j(x)h(x)\dx\Big|\\
&=\Big|\int t^{2M}(r^2_{B_j}L)^kt^2L\Phi(t\sqrt{L})(a_j(\cdot,t))(x)h(x)\dtt \dx\Big|\\
&\leq Cr^{2M}_{B_j}\Big(\int_{\XL}|a_j(x,t)|^2
\f{w(B(x,t))}{V(x,t)}\f{V(x,t)}{w(B(x,t))}\dtt \dx\Big)^{1/2}\\
&~~~\times\Big(\int_{\XL} |(t^2L)^{k+1}\Phi(t\sqrt{L})(x)h(x)|^2\dtt
\dx\Big)^{1/2}.
\end{aligned}
\end{equation*}
Since $(x,t)\in \widehat{B_j}$, the ball $B(x,t)\subset B_j$. This together with Lemma \ref{weightedlemma2} gives
 \begin{equation*}
\begin{aligned}
\Big|\int (r^2_{B_j}L)^kb_j(x)h(x)\dx\Big|
&\leq  Cr^{2M}_{B_j}\f{V(B_j)}{w(B_j)}\|a\|_{L^2(W)}\|h\|_{L^2(B_j)}\\
&\leq Cr^{2M}_{B_j}\f{V(B_j)^{1/2}}{w(B)^{1/2}}w(B_j)^{1/2-1/p}\|h\|_{L^2(B_j)}\\
&=Cr^{2M}_{B_j}V(B_j)^{1/2}w(B_j)^{-1/p}\|h\|_{L^2(B_j)}.
\end{aligned}
\end{equation*}
This implies $\pi_{\Psi,L}a_j$ is a multiple of an $(M,p,w)$ atom with the harmless constant
for each $j$.

To complete the proof, it is suffice to show that for any $(M,p,w)$
atom $a$ in $\HL$ we have
\begin{equation}\label{atomisinHardyspace}
\|S_La\|_{L^p_w(X)}\leq C.
\end{equation}
Indeed, we write
\begin{equation*}
\begin{aligned}
\|S_La\|^p_{L^p_w(X)}&=\int_\RR|S_La(x)|^pw(x)\dx\\
&=\sum_{j=0}^\vc \int_{S_j(B)} |S_La(x)|^pw(x)\dx\\
&=\sum_{j=0}^\vc I_j.
\end{aligned}
\end{equation*}
For $j=0,1,2$ we have, by the H\"older inequality, $L^2$-boundedness of
$S_L$, and $w\in A_1\cap RH_{\f{2}{2-p}}$,
\begin{equation*}
\begin{aligned}
I_j&\leq \|S_L a\|_{L^2(S_j(B))}^p\Big(\int_{S_j(B)}w^{\f{2}{2-p}}\Big)^{\f{2-p}{2}}\\
&\leq C\|a\|^p_{L^2}V(2^jB)^{-p/2}w(2^jB)\\
&\leq CV(B)^{p/2}w(B)^{-1}V(2^jB)^{-p/2}w(2^jB).
\end{aligned}
\end{equation*}
Since $w\in A_1$, using Lemma \ref{weightedlemma2}, we have
$$
I_j\leq CV(B)^{p/2 -1}V(2^jB)^{1-p/2}\leq C \ \text{for $j=0,1,2$}.
$$
 For $j\geq 3$, we have
$$
I_j\leq \|S_L a\|_{L^2(S_j(B))}^p\Big(\int_{S_j(B)}w^{\f{2}{2-p}}\Big)^{\f{2-p}{2}}\leq C\|S_L a\|_{L^2}^pV(2^jB)^{-p/2}w(2^jB)
$$
To estimate $\|S_L a\|_{L^2(S_j(B))}^p$, we write
\begin{equation*}
\begin{aligned}
\|S_L a&\|_{L^2(S_j(B))}^2\\
&=\int_{S_j(B)}\Big(\int_{0}^{\f{d(x,x_B)}{4}}+\int^\vc_{\f{d(x,x_B)}{4}}\Big)\int_{d(x,y)<t}|t^2Le^{-t^2L}a(y)|^2\dyt \dx\\
&=\int_{S_j(B)}\Big(\int_{0}^{\f{d(x,x_B)}{4}}+\int^\vc_{\f{d(x,x_B)}{4}}\Big)\int_{d(x,y)<t}|(t^2L)^{M+1}e^{-t^2L}b(y)|^2\f{\dy}{V(x,y)}\f{dt}{t^{4M+1}} \dx\\
&=J_1+J_2,
\end{aligned}
\end{equation*}
where $a=L^Mb$.\\
Setting $F_j(B):=\{y:d(x,y)<\f{d(x,x_B)}{4} \ \text{for some $x\in
S_j(B)$}\}$, then $d(B, F_j(B))\geq 2^{j-2}r_B$. By the fact that
$\int_{d(x,y)<t}\f{1}{V(x,t)}\dx <C$, we have
\begin{equation*}
\begin{aligned}
J_1&\leq \int_{0}^{2^{j}r_B}\int_{F_j(B)}|(t^2L)^{M+1}e^{-t^2L}b(y)|^2\dy\f{dt}{t^{4M+1}}\\
&\leq \|b\|^2_{L^2}\int_{0}^{2^{j}r_B}\exp\Big(-\f{d^2(B, F_j(B))}{ct^2}\Big)\f{dt}{t^{4M+1}}\\
&\leq r_B^{4M}V(B)w(B)^{-2/p}\int_{0}^{2^{j}r_B}\Big(\f{t}{2^jr_B}\Big)^{4M+1}\f{dt}{t^{4M+1}}\\
&\leq 2^{-4jM}V(B)w(B)^{-2/p}.
\end{aligned}
\end{equation*}
For the term $J_2$ we have
\begin{equation*}
\begin{aligned}
J_2&\leq \int_{2^{j-1}r_B}^\vc|(t^2L)^{M+1}e^{-t^2L}b(y)|^2\dy\f{dt}{t^{4M+1}}\\
&\leq \|b\|^2_{L^2}\int_{2^{j-1}r_B}^\vc\f{dt}{t^{4M+1}}\\
&\leq 2^{-4jM}V(B)w(B)^{-2/p}.
\end{aligned}
\end{equation*}
Combining above estimates of $J_1$ and $J_2$, using $w\in A_1$ we
have
\begin{equation*}
\begin{aligned}
I_j&\leq C2^{-2jpM}V(2^jB)^{-p/2}w(2^jB)V(B)^{p/2}w(B)^{-1}\\
&\leq C2^{-2jpM}V(2^jB)^{-p/2}V(2^jB)V(B)^{p/2}V(B)^{-1}\\
&\leq C2^{-2jMp}2^{-jnp/2}2^{jn}=C2^{-jp(2M+\f{n}{2}-\f{n}{p})}.
\end{aligned}
\end{equation*}
Since $M>\f{n(2-p)}{4p}$, $\|S_La\|_{L^p_w(X)}\leq C$.

\begin{thm}\label{thm1-atomicdecompostion} Let $L$ satisfy ($H1$) and ($H2$).
Let $M \geq \f{n(2-p)}{4p}$ and $w \in A_1$. If $f \in \HL\cap
H^2(X)$, then there exists a family of $(M, p, w)$-atoms
$\{a_j\}_{j=0}^\vc$ and a sequence of numbers
$\{\lambda_j\}_{j=0}^\vc$ such that $f$ can be represented in the
form $f = \sum_{j=0}^\vc\lambda_ja_j$, and the sum converges in the
sense of $L^2(X)$-norm. Moreover,
$$
\sum_{j=0}^\vc|\lambda_j|^p\leq C\|f\|^p_{\HL}.
$$
\end{thm}
\emph{Proof:} Let $f\in f \in \HL\cap H^2(X)$. Then
$t^2Le^{-t^2L}f\in T^p_w(X)\cap T^2(X)$. Then we obtain
$$
f=\pi_{\Psi, L}(t^2Le^{-t^2L}f).
$$
By using the argument as in Proposition
\ref{Pro1-mappingfrom tentspace to Hardyspace} we complete the
proof.\\

\begin{thm}\label{thm2-atomicdecom}
Given $w \in A_1\cap RH_{\f{2}{2-p}}$ and $M>\f{n(2-p)}{4p}$. Let $f =
\sum_{j=0}^\vc\lambda_ja_j$, where $\{\lambda_j\}_{j=0}^\vc\in l^p$,
$a_j$'s are $(M, p, w)$-atoms, and the sum converges in $L^2(X)$.
Then $f\in \HL$ and
    $$
    \Big\|\sum_{j=0}^\vc\lambda_ja_j \Big\|^p_{\HL}\leq C\sum_{j=0}^\vc|\lambda_j|^p.
    $$
\end{thm}
\emph{Proof:} Since $S_L$ is bounded on $L^2(X)$, $S_Lf\leq |\lambda_j|\sum_{j}S_La_j$. Therefore, to show $f\in \HL$, it is sufficient to claim that there exists a constant $C>0$ so that
$$
\|S_La\|_{L^p_w(X)}\leq C
$$
for all $(M,p,w)$ atoms $a$.

This have been proved in (\ref{atomisinHardyspace}) and hence we complete the proof.\\

As a consequence of Theorems \ref{thm1-atomicdecompostion} and \ref{thm2-atomicdecom}, we conclude that
\begin{cor}
Let $0<p\leq 1$, $w\in A_1\cap RH_{\f{2}{2-p}}$ and $M>\f{n(2-p)}{4p}$. Then the spaces $\HL$ and $H^p_{L,M,w,at}(X)$ coincide and their norms are equivalent.
\end{cor}

\subsection{Molecular characterization of weighted Hardy spaces $\HL$}
The weighted Hardy spaces $\HL$ can be also characterized in
molecular decomposition. More specifically, we have the following result.
\begin{thm} Let $L$ satisfy ($H1$) and ($H2$).
\begin{enumerate}[(i)]
\item Let $M \geq 1$ and $w \in A_1\cap RH_{\f{2}{2-p}}$. If $f \in \HL\cap H^2(X)$,
then there exists a family of $(M, p, w, \epsilon)$-molecules
$\{m_j\}_{j=0}^\vc$ and a sequence of numbers
$\{\lambda_j\}_{j=0}^\vc$ such that $f$ can be represented in the
form $f = \sum_{j=0}^\vc\lambda_jm_j$, and the sum converges in the
sense of $L^2(X)$-norm. Moreover,
$$
\sum_{j=0}^\vc|\lambda_j|^p\leq C\|f\|^p_{\HL}.
$$
\item Conversely, given $w \in A_1\cap RH_{2}$ and $M>\f{n(2-p)}{4p}$. Let $f =
\sum_{j=0}^\vc\lambda_ja_j$, where $\{\lambda_j\}_{j=0}^\vc\in l^p$,
$m_j$'s are $(M, p, w, \epsilon)$-molecules, and the sum converges
in $L^2(X)$. Then $f\in \HL\cap L^2(X)$ and
    $$
    \Big\|\sum_{j=0}^\vc\lambda_jm_j \Big\|^p_{\HL}\leq C\sum_{j=0}^\vc|\lambda_j|^p.
    $$
    \end{enumerate}
\end{thm}

\emph{Proof:} (i) The proof of (i) is a direct consequence of
Theorem \ref{thm1-atomicdecompostion}, since an  $(M, p, w)$-atom is
also an  $(M, p, w, \epsilon)$-molecule for all $\epsilon >0$.

\medskip

(ii) The proof of (ii) is  similar to that of Theorem
\ref{thm2-atomicdecom}. The main difference is that  the support of
$(M, p, w, \epsilon)$-molecule is not the ball $B$. However, we can
overcome this difficulty by decomposing $X$ into annuli associated
with the ball $B$, then using the same argument as in Theorem
\ref{thm2-atomicdecom} to get (ii). We omit the details here.

\section{Boundedness of singular integrals with non-smooth kernels}

In this section, we study the boudnedness of some singular integrals such as the generalized Riesz transforms and spectral multipliers on the weighted Hardy spaces. Before coming to details, we need the following technical lemma.
\begin{lem}\label{lem1-characterizationofHardyspace}
Suppose that $T$ is a linear (resp. nonnegative sublinear) operator
which maps $L^2(X)$ continuously into $L^{2,\vc}(X)$. Let $0<p\leq 1$ and $w\in A_1\cap RH_{\f{2}{2-p}}$.
\begin{enumerate}[(i)]
\item If there exists  a constant
$C$ such that
$$
||Ta||_{L_w^{p}(X)}\leq C
$$
for all $(M,p,w)$-atoms $a$, then $T$ extends to a bounded linear
(resp. sublinear) operator from $H^p_{L,w}(X)$ to $L^{p}(w)$.
\item If there exists a constant
$C$ such that
$$
||Ta||_{\HL}\leq C
$$
for all $(M,p,w)$-atoms $a$, then $T$ extends to a bounded linear
(resp. sublinear) operator on $H^p_{L,w}(X)$.
\end{enumerate}
\end{lem}
The proof of this lemma is similar to that in \cite[Lemma 4.1]{HM} and hence we omit details here.

\subsection{Generalized Riesz transforms}

Assume that  $L$ satisfies (H1) and (H2). Also assume that $D$ is a densely defined
linear operator on $L^2(X)$ which possesses the following
properties:
\begin{enumerate}[(i)]
\item $DL^{-1/2}$ is bounded on $L^2$;
\item the family operators $\{\sqrt{t}De^{-tL}\}_{t>0}$ satisfies the
Davies-Gaffney estimate (H2).
\end{enumerate}
Examples of such an  operator $D$  include gradient operator in
the Euclidean space and the Riemannian gradient on complete
Riemannian manifolds, see for example \cite{A,ACDH,CD1}.

\begin{thm}\label{thm11}
For any $f\in \HL$ for $0<p\leq 1$ and $w\in A_1\cap RH_{\f{2}{2-p}}$,
$$
||DL^{-1/2} (f)||_{L^p_w(X)}\leq C||f||_{\HL}.
$$
\end{thm}
Before giving the proof of Theorem \ref{thm11}, we state the
following lemma.
\begin{lem}\label{lem3}
For every $M\in \mathbb{N}$, all closed sets E, F in $X$ with
$d(E,F)>0$ and every $f\in L^2(X)$ supported in $E$, one has
\begin{equation}\label{eq2}
||DL^{-1/2}(I-e^{-tL})^Mf||_{L^2(F)}\leq
C\Big(\f{t}{d(E,F)^2}\Big)^M||f||_{L^2(E)}, \ \forall t>0,
\end{equation}
and
\begin{equation}\label{eq3}
||DL^{-1/2}(tLe^{-tL})^Mf||_{L^2(F)}\leq
C\Big(\f{t}{d(E,F)^2}\Big)^M||f||_{L^2(E)}, \ \forall t>0.
\end{equation}
\end{lem}
\emph{Proof:} The proof is similar to that of Lemma
2.2 in \cite{HMa} and we omit it here.\\

\noindent We now prove Theorem \ref{thm11}.

\medskip

\noindent \emph{Proof of Theorem \ref{thm11}:}
Due to Lemma \ref{lem1-characterizationofHardyspace}, we need only to claim that there exists a constant $C>0$ so that
$$
\|DL^{-1/2}a\|_{L^p_w(X)}\leq C
$$
for all $(M,p,w)$-atoms associated to a ball $B$ with $M>\f{n(2-p)}{4p}$.

Indeed, let $b\in L^2(X)$ so that $a=L^Mb$. Setting $T=DL^{-1/2}$, we write
\begin{equation*}
\begin{aligned}
\|Ta\|^p_{L^p_w(X)}=\int_{X}& |Ta(x)|^pw(x)d\mu(x)\\
&\leq
\int_{X}|T((I-e^{r_B^2L})^Ma)(x)|^pw(x)d\mu(x)\\
&~~~~~~~~~+\int_{X}|T([I-(I-e^{r_B^2L})^M]a)(x)|^pw(x)d\mu(x)\\
&= I+II.
\end{aligned}
\end{equation*}
We estimate the term $I$ first. By the H\"older inequality, we obtain
\begin{equation*}
\begin{aligned}
I&\leq
\sum_{k=0}^{\vc}\int_{S_k(B)}|T((I-e^{r_B^2L})^Ma)(x)|^pw(x) d\mu(x)\\
&\leq \sum_{k=0}^{\vc}
\|T((I-e^{r_B^2L})^Ma)\|^p_{L^2(S_k(B))}\Big(\int_{S_k(B)}w^{\f{2}{2-p}}d\mu\Big)^{\f{2-p}{2}}\\
&\leq \sum_{k=0}^{\vc}
\|T((I-e^{r_B^2L})^Ma)\|^p_{L^2(S_k(B))}\int_{S_k(B)}V(2^kB)^{-p/2}w(2^kB):=\sum_{k=0}^\vc
I_k.
\end{aligned}
\end{equation*}
For $k=0,1,2,$ one has
$$
I_k\leq C\|a\|^p_{L^2(B)}V(2^kB)^{-p/2}w(2^kB)\leq
CV(B)^{p/2}w(B)^{-1}V(2^kB)^{-p/2}w(2^kB).
$$
This together with Lemma \ref{weightedlemma2} implies, for $k=0,1,2$,
$$
I_k\leq CV(B)^{p/2-1}V(2^kB)^{1-p/2}\leq C .
$$
For $k\geq 3$
\begin{equation*}
\begin{aligned}
||T((I-e^{r_B^2L})^M&a)||_{L^2(S_k(B))}\leq
C2^{-2Mk}||a||_{L^2(B)}\leq C2^{-2Mk}V(B)^{1/2}w(B)^{-1/p}.
\end{aligned}
\end{equation*}
Therefore,
\begin{equation*}
\begin{aligned}
I_k&\leq C2^{-2Mpk}V(B)^{p/2}w(B)^{-1}V(2^kB)^{-p/2}w(2^kB)\\
&\leq C2^{-2Mpk}V(B)^{p/2}V(B)^{-1}V(2^kB)^{-p/2}V(2^kB)\\
&\leq C2^{-2Mpk}V(B)^{p/2-1}V(2^kB)^{1-p/2}\\
&\leq C2^{-2Mpk+kn(1-p/2)}.
\end{aligned}
\end{equation*}
Due to $M>\f{n(2-p)}{4p}$, we obtain $\sum_{k=0}^\vc I_k\leq C$.

\medskip

\noindent For the term $II$, the same argument above gives
\begin{equation*}
\begin{aligned}
II&\leq
\sum_{k=0}^{\vc}\int_{S_k(B)}|T([I-(I-e^{r_B^2L})^M]a)(x)|^pw(x)d\mu(x)\\
&\leq
\sum_{k=0}^{\vc}\|T([I-(I-e^{r_B^2L})^M]a)\|_{L^2(S_k(B))}V(2^kB)^{-p/2}w(2^kB)\\
&\leq \sum_{k=0}^{\vc}II_k.
\end{aligned}
\end{equation*}
Next we have
$$
I-(I-e^{r_B^2L})^M = \sum_{l=1}^{M}c_le^{-lr_B^2L},
$$
where $c_l=(-1)^{l+1}\f{M!}{(M-l)!l!}$. Therefore,
\begin{equation*}
\begin{aligned}
II_k &\leq C\sup_{1\leq l\leq
M}\|Te^{-lr_B^2L}a\|_{L^2(S_k(B))}V(2^kB)^{-p/2}w(2^kB)\\
&\leq C\sup_{1\leq l\leq
M}\Big|\Big|T\Big(\f{l}{M}r_B^2Le^{-\f{l}{M}r_B^2L}\Big)^M(r_B^{-2}L^{-1})^Ma\Big|\Big|_{L^2(S_k(B))}V(2^kB)^{-p/2}w(2^kB).
\end{aligned}
\end{equation*}
At this point, by the same argument as in the estimate $I_k$, we
also obtain that $II<C$. This therefore completes our proof.\\

\subsection{Boundedness of Riesz transforms associated with  magnetic Schr\"odinger operators}

\subsubsection{  Magnetic Schr\"odinger operators and heat kernel estimates}

Consider magnetic Schr\"odinger operators in general setting as in
\cite{DOY}. Let the real vector potential ${\vec a}=(a_1, \cdots,
a_n)$ satisfy
\begin{eqnarray}
a_k\in L^2_{\rm loc}({\mathbb R}^n), \ \ \ \ \forall k=1, \cdots, n,
\ \ \label{e1.1}
\end{eqnarray}
and an electric potential $V$ with
\begin{eqnarray}
0\leq V\in L^1_{\rm loc}({\mathbb R}^n). \label{e1.2}
\end{eqnarray}

\noindent Let $L_k={\partial/\partial x_k}-i a_k $. We   define the
form $Q$ by
\begin{eqnarray*}
Q(f,g)=\sum_{k=1}^n\int_{{\mathbb R}^n}L_kf
   {\overline {L_kg}} \  \! dx + \int_{{\mathbb R}^n}V f
   {\overline {g}} \ \! dx
  \end{eqnarray*}

  \noindent
with domain $D(Q)=\mathcal{Q}\times \mathcal{Q} $ here
$$ \mathcal{Q}=\{f\in L^2({\mathbb R}^n), L_kf\in L^2({\mathbb R}^n) {\rm \  for}\
k=1,\cdots, n {\rm  \ and} \ \sqrt{V} f\in L^2({\mathbb R}^n)\}.
$$

 \noindent
 It is well known that this symmetric form is closed and this form coincides with the minimal closure
of the form given by the same expression but defined on
$C^{\infty}_0({\mathbb R}^n)$ (the space of $C^{\infty}$ functions
with compact supports). See, for example \cite{Si}.

Let us denote by $A$ the non-negative self-adjoint operator associated with $Q$.
The domain of $A$ is given by
$$
{\mathcal D}(A)=\Big\{f\in {\mathcal D}(Q), \exists g\in L^2({\Bbb
R}^n) {\rm \ such \ that\ } Q(f,\varphi)=\int_{{\mathbb R}^n} g{\bar
\varphi} dx, \ \ \forall\varphi\in {\mathcal D}(Q)\Big\},
$$
and $A$ is given by the expression
\begin{eqnarray}
Af=\sum_{k=1}^nL_k^{\ast} L_kf+Vf. \label{e1.3}
\end{eqnarray}

\noindent Formally, we write  $A=-(\nabla-i{\vec a})\cdot
(\nabla-i{\vec a})+V$. The operator $A$ generates a semigroup $e^{-tA}$ which possesses
a Gaussian kernel bound. Indeed, by the well known diamagnetic inequality
(see, Theorem 2.3 of \cite{Si} and \cite{CFKS} for instance) we have
the pointwise inequality
\begin{eqnarray*}
\big|e^{-tA}f(x)\big|\leq e^{t\triangle}\big(|f|\big)(x) \ \ \
\forall t > 0, \ \ f\in L^2({\mathbb R}^n).
\end{eqnarray*}
\noindent This inequality implies in particular that the semigroup
$e^{-tA}$ maps $L^1({\Bbb R}^n)$ into $L^{\infty}({\Bbb R}^n)$ and
that the kernel $p_t(x,y)$ of $e^{-tA}$ satisfies
\begin{eqnarray}
  \big|p_t(x,y)\big|\leq  (4\pi t)^{-{n\over 2}}
 \exp\Big(-\frac{|x-y|^2}{4t}\Big)
\label{e2.5}
\end{eqnarray}

\noindent
  for all $t>0$ and  almost all $x,y\in {\mathbb R}^n$.\\

For  $k=1, \cdots, n$, the operators
$L_kA^{-1/2}$ are called the Riesz transforms associated with $A.$
It is easy to check that
\begin{eqnarray}
 \|L_kf \|_{L^2({\mathbb R}^n)}\leq  \|A^{1/2} f \|_{L^2({\mathbb R}^n)},
 \ \ \ \ \ \  \  \forall f\in {\mathcal
D}(Q)={\mathcal D}(A^{1/2}) \label{e1.4}
\end{eqnarray}

\noindent for any  $k=1, \cdots, n$, and hence the operators
$L_kA^{-1/2}$ are bounded on $L^2({\mathbb R}^n)$. Note that this is
also true for $V^{1/2}A^{-1/2}$. Moreover, it was recently proved in
Theorem 1.1 of  \cite{DOY}   that for each $k=1, \cdots, n$,  the
Riesz transforms $L_k A^{-1/2}$ and $V^{1/2} A^{-1/2}$ are   bounded
on $L^p({\mathbb R}^n)$ for all $1<p\leq 2$, i.e., there exists a
constant $C_p>0$ such that
\begin{eqnarray}
\label{ev}\hspace{1cm}
    \big\|V^{1/2} A^{-1/2}f\big\|_{L^p({\mathbb R}^n)}+
    \sum_{k=1}^n  \big\|L_k A^{-1/2}f\big\|_{L^p({\mathbb R}^n)} \leq C_p\|f\|_{L^p({\mathbb R}^n)}, \ \ \
\end{eqnarray}
for $1<p\leq 2$.

\noindent
The $L^p$-boundedness of Riesz transforms for the range $p>2$ can be obtained
if one imposes certain additional regularity conditions on the potential $V$, see for example \cite{AB}.\\

Let us recall that
\begin{enumerate}[(i)]

\item In \cite{DOY}, the boundedness of the Riesz transforms
$L_kA^{-1/2}$ and $V^{1/2}A^{-1/2}$ was proved for
$L^{p}(\mathbb{R}^n)$ spaces with $1<p<2$;
\item Recently, \cite{A} extended the results in \cite{DOY} to
weighted weak type $L^{1,\infty}$ estimates and weighted $L^p$
estimates with an appropriate range of $p$ (depending on the weight).
\item For $p\in (0,1]$, it was proved that the Riesz transforms
$L_kA^{-1/2}$ and $V^{1/2}A^{-1/2}$ are bounded from $H_A^p$ to
$L^p$, where $H^p_A$ denotes the unweighted Hardy space associated to operator
$A$, see \cite{AD}.

\end{enumerate}

It is natural to ask the question of the boundedness of
the Riesz transforms $L_kA^{-1/2}$ and $V^{1/2}A^{-1/2}$ on the Hardy spaces with weights on the
range $0<p\leq 1$. Our aim in this section is to give a positive answer and
establish the weighted estimates for Riesz transforms
 $L_kA^{-1/2}$ and $V^{1/2}A^{-1/2}$ for $0<p\leq 1$. The main result
of this section  is the following theorem.
\begin{thm}\label{thm4}
If $w\in A_1\cap RH_{\f{2}{2-p}}$ then the Riesz transforms $\LA$ and $\VA$ are bounded from $H^p_{A,w}(\mathbb{R}^n)$ to $L^p_w(\mathbb{R}^n)$ for all $0<p\leq 1$.
\end{thm}

 \subsubsection{Proof of Theorem \ref{thm4}}
The following propositions give  the estimates for the operators $\LA$
and $\VA$ which will be useful for the proof of Theorem \ref{thm4}.
\begin{prop}\label{prop2.1} For all $m\geq 1$ there exists $C>0$ such that
for all balls $B$ and all $f$ with support in $B$
\begin{equation}\label{eq4.1}
\Big(\int_{S_j(B)}|\VA(I-e^{-r_B^2A})^m f|^2\
dx\Big)^{\frac{1}{2}}\leq C\ 2^{-2j(m-1)}\Big(\int_B|f|^2\Big)^{1/2}
\end{equation}
and for all $k$,
\begin{equation}\label{eq4.2}
\Big(\int_{S_j(B)}|\LA(I-e^{-r_B^2A})^m f|^2\
dx\Big)^{\frac{1}{2}}\leq C\ 2^{-2j(m-1)}\Big(\int_B|f|^2\Big)^{1/2} .
\end{equation}
\end{prop}

\begin{prop}\label{prop2.2} For all $m\geq 1$ and $t >  0$ there exists $C>0$ such that
for all balls $B$, $\f{1}{m}r_B^2\leq t\leq r_B^2$ and all $f$ with
support in $B$
\begin{equation}\label{eq5.1}
\Big(\int_{S_j(B)}|\VA(tAe^{-tA})^m f|^2\
dx\Big)^{\frac{1}{2}}\leq C2^{-2jm}\Big(\int_B|f|^2\Big)^{1/2}
\end{equation}
and for all $k$,
\begin{equation}\label{eq5.2}
\Big(\int_{S_j(B)}|\LA (tAe^{-tA})^m f|^2\
dx\Big)^{\frac{1}{2}}\leq C2^{-2jm}\Big(\int_B|f|^2\Big)^{1/2} .
\end{equation}
\end{prop}

The proofs of Propositions \ref{prop2.1} and \ref{prop2.2} are standard and rely on the heat kernel estimates and   we refer the reader to \cite{AD}
for the details of the proofs.
From Propositions \ref{prop2.1} and \ref{prop2.2}, the following estimates hold.
\begin{prop}\label{prop2.3} There exists $C>0$ such that for any $(M,p,w)$ atoms $a$ associated
to the ball $B$ we have
\begin{equation}\label{e4.1}
\Big(\int_{S_j(B)}|\VA a|^2\
dx\Big)^{\frac{1}{2}}\leq C\ 2^{-2j(M-1)}V(B)^{1/2}w(B)^{-1/p}
\end{equation}
and for all $k$,
\begin{equation}\label{e4.2}
\Big(\int_{S_j(B)}|\LA a|^2\
dx\Big)^{\frac{1}{2}}\leq C\ 2^{-2j(M-1)}V(B)^{1/2}w(B)^{-1/p}.
\end{equation}
\end{prop}

\emph{Proof:} Set $T=\VA$ (or, $T=\LA$). We have, for each $(M,p,w)$-atom $a=A^Mb$ associated to the ball $B$,
$$Ta=T(I-e^{-r_B^2A})^Ma+ T(I-(I-e^{-r_B^2A})^M)a .$$
Observe that
$$ I-(I-e^{-r_B^2A})^M =
\sum_{k=1}^{M}c_ke^{-kr_B^2A},
$$
where $c_k=(-1)^{k+1}\f{M!}{(M-k)!k!}$. Therefore,
\begin{equation*}
T[I-(I-e^{-r_B^2A})^M]a =
\sum_{k=1}^{M}a_kT\Big(r_B^2\f{k}{M}A
e^{-\f{k}{M}r_B^2A}\Big)^M (r_B^{-2M}b),
\end{equation*}
where $a_k=c_k\Big(\f{k}{M}\Big)^M$.\\
At this stage, applying Propositions \ref{prop2.1} and  \ref{prop2.2}, we obtain (\ref{e4.1}) and (\ref{e4.2}).\\

\emph{Proof of Theorem \ref{thm4}:} By Lemma \ref{lem1-characterizationofHardyspace},
it suffices to show that for any $(M,p, w)$-atom $a$ associated to the ball $B$ with
$M>1+f{n(2-p)}{4p}$,   we have $||T a||_{L_w^p}\leq C$ where $T=\VA$ or $T=\LA$.\\
Indeed, we have
\begin{equation*}
\begin{aligned}
||T a||_{L_w^p}^p=\int_X|T a|^pw(x) dx &= \sum_{j=0}^\infty \int_{S_j(B)}|T a|^pw(x) dx \\
&= \sum_{j=0}^\infty K_j.
\end{aligned}
\end{equation*}
By the H\"older inequality and (\ref{e4.1}), (\ref{e4.2}), one has, for each $j$,
\begin{equation*}
\begin{aligned}
K_j&\leq ||Ta||^p_{L^2(S_j(B))}\Big(\int_{S_j(B)}\Big(\int_{S_j(B)} w^{\f{2}{2-p}}\Big)^{\f{2-p}{2}}\\
& \leq C 2^{-2jp(M-1)}V(B)^{p/2}w(B)^{-1}V(2^jB)^{-p/2}w(2^jB)\\
& \leq C 2^{-jp[2(M-1)+\f{n}{2}-\f{n}{p}]}.\\
\end{aligned}
\end{equation*}
This, together with $M>1+\f{n}{p}-\f{n}{2}$, gives
$$
||T a||_{L_w^p}^p\leq C.
$$
 The proof is complete.

\subsection{Spectral multiplier theorem on $H_{L,w}^p(X)$}
In this section, by the kernel $K_T(x, y)$ associated to a $L^2$-bounded linear operator $T$ we
mean
$$
Tf(x) =\int_X K_T(x, y)f(y)d\mu(y)
$$
where $K_T(x, y)$ is a measurable function and the formula above holds for each continuous
function $f$ with bounded support and for almost all $x$ not in the support of $f$.

Our main results are the following two theorems.
\begin{thm}\label{thm1}
Let $L$ be an operator satisfying $(H1)$ and $(H3)$, and $w\in
A_1\cap RH_{\f{2}{2-p}}$. Suppose that $s>\f{n(2-p)}{2p}$ and for any $R>0$ and all
Borel functions $F$ such that {\rm supp}\,$F\subset [0,R]$,
\begin{equation}\label{eq3.1}
\int_X |K_{F(\sqrt{L})}(x,y)|^2\dx\leq
\f{C}{V(y,R^{-1})}\|\delta_RF\|^2_{L^q}
\end{equation}
for some $q\in [2,\infty]$. Then for any   function $F$ such
that $\sup_{t>0}\|\eta \delta_tF\|_{W^q_s}<\infty$, where $ \delta_t F(\lambda)=F(t\lambda)$,
   $\| F \|_{W^q_s}=\|(I-d^2/d x^2)^{s/2}F\|_{L^q}$, the operator
$F(L)$ is bounded on $\HL$.
\end{thm}

Note that (\ref{eq3.1}) always holds for $q=\vc$, see \cite{DOS}.
For further discussions concerning condition (\ref{eq3.1}), we refer the reader
to \cite{DOS}.\\

As a preamble to the proof of Theorem ~\ref{thm1}, we record a
useful auxiliary result which is taken from \cite[Lemma 4.3]{DOS}.

\begin{lem}\label{lem1}
Suppose that $L$ satisfies (\ref{eq3.1}) for some $q\in [2, \infty]$,
 $R>0$ and $ s>0$. Then for any $\epsilon>0$, there exists a
constant $C=C(s, \epsilon)$ such that
\begin{eqnarray}
\int_X \big|K_{F(\sqrt{L})}(x,y)\big|^2 \big(1+Rd(x,y)\big)^{s}
d\mu(x)\leq \f{C}{ V(y, R^{-1})}
 \|\delta_{R} F\|^2_{W^q_{\f{s}{2} +\epsilon}}
\end{eqnarray}
for all Borel functions $F$ such that {\rm supp}\,$F\subseteq [R/4,
R].$
\end{lem}

{\it Proof of Theorem \ref{thm1}:} Since condition $\sup_{t>0}\|\eta
\delta_tF\|_{W^q_s}<\infty$ is invariant under the change of
variable $\lambda \mapsto \lambda^s$ and independent on the choice
of $\eta$, due to Lemma \ref{lem1-characterizationofHardyspace} it
suffices to show that there exists $\epsilon>0$ such that for any
$(2M,p,w)$-atom $a=L^{2M}b$ in $\HL$ the function
$$
\widetilde{a}=F(\sqrt{L})a=L^M(F(\sqrt{L})L^Ma)
$$
is a multiple of $(M,p,w,\epsilon)$-molecule for $M>\f{n(2-p)}{4p}$ and some $\epsilon>0$. By definition, it suffices to prove that for all $l=0,1,\ldots, M$,
\begin{equation}\label{eq1-spectralmultipliers}
\|(r_B^2L)^l(F(\sqrt{L})L^Ma)\|_{L^2(S_k(B))}\leq C2^{-k\epsilon}r_B^{2M}V(B)^{1/2}w(B)^{-1/p}
\end{equation}
for all $k=0,1,2, \ldots$.\\

For $k=0,1,2$, using $L^2$-boundedness of $F(\sqrt{L})$, we have
$$
\|(r_B^2L)^l(F(\sqrt{L})L^Ma)\|_{L^2(S_k(B))}\leq C\|(r_B^2L)^lL^Ma)\|\leq C2^{-k\epsilon}r_B^{2M}V(B)^{1/2}w(B)^{-1/p}.
$$

Now we need only to verify (\ref{eq1-spectralmultipliers}) for $k\geq 3$. By standard argument, fix a function $\phi \in
C_c^\infty(\f{1}{4},1)$ such that
$$
\sum_{j\in \mathbb{Z}}\phi(2^{-j}\lambda)=1 \ \text{for} \
\lambda>0.
$$
Then, for $0\leq l\leq M$, one has
\begin{equation}\label{eq1}
\begin{aligned}
(r_B^2L)^l\widetilde{b}&=r_B^{2l}\sum_{j\geq j_0}\phi(2^{-j}\sqrt{L})F(\sqrt{L})L^{l+M}b\\
&~~~+r_B^{2l}\sum_{j< j_0}\phi(2^{-j}\sqrt{L})L^M F(\sqrt{L})L^{l}b\\
&=r_B^{2l}\sum_{j\geq j_0}\phi(2^{-j}\sqrt{L})
F(\sqrt{L})b_1+r_B^{2l}\sum_{j< j_0}\phi(2^{-j}\sqrt{L})L^M
F(\sqrt{L})b_2,
\end{aligned}
\end{equation}
where $\widetilde{b}=L^Mb$ and $j_0=-\log_2r_B$.\\
It is easy to see that
$$
\|b_1\|_{L^2}\leq r_B^{2M-2l}V(B)^{\f{1}{2}}w(B)^{-1/p} \ \text{and}
\ \|b_2\|_{L^2}\leq r_B^{4M-2l}V(B)^{\f{1}{2}}w(B)^{-1/p}.
$$
Setting
$$
F_j(\lambda)=
\begin{cases}
F(\lambda)\phi(2^{-j}\lambda) \ & \text{for} \ \ \ j\geq j_0\\
F(\lambda)(2^{-j}\lambda)^{2M}\phi(2^{-j}\lambda) \ & \text{for} \ \
\  j< j_0,
\end{cases}
$$
then we can rewrite (\ref{eq1}) as follows
\begin{equation}\label{eq2}
\begin{aligned}
(r_B^2L)^l\widetilde{b}&=r_B^{2l}\sum_{j\geq
j_0}F_j(\sqrt{L})b_1+r_B^{2l}2^{2jM}\sum_{j< j_0}F_j(\sqrt{L})b_2.
\end{aligned}
\end{equation}
Since (\ref{eq2}) converges in $L^2(X)$, we have for any $k\geq 3$
\begin{equation*}
\begin{aligned}
\|(r_B^2L)^l\widetilde{b}\|_{L^2(S_k(B))}\leq  r_B^{2l}\sum_{j\geq
j_0}\|F_j(\sqrt{L})b_1\|_{L^2(S_k(B))}+r_B^{2l}2^{2jM}\sum_{j<
j_0}\|F_j(\sqrt{L})b_2\|_{L^2(S_k(B))}.
\end{aligned}
\end{equation*}

Let us estimate $\|F_j(\sqrt{L})b_1\|_{L^2(S_k(B))}$ for $j\geq j_0$
first. Since supp\,$F_j\subset [R/4,R]$ with $R=2^j$, applying Lemma
\ref{lem1} and the Minkowski inequality, we have, for
$s>s'>\f{n(2-p)}{2p}\geq \f{n}{2}$ and $k\geq 3$,
\begin{equation}\label{eq3}
\begin{aligned}
\|F_j(\sqrt{L})&b_1\|_{L^2(S_k(B))}\\
&\leq \Big\|\int_{B}K_{F_j(\sqrt{L})}(x,y)b_1(y)\dy\Big\|_{L^2(S_k(B))}\\
&\leq \|b_1\|_{L^1}\sup_{y\in B}\Big(\int_{S_k(B)}|K_{F_j(\sqrt{L})}(x,y)|^2\dx\Big)^{1/2}\\
&\leq \|b_1\|_{L^2}V(B)^{\f{1}{2}}\sup_{y\in B}\Big(\int_{S_k(B)}|K_{F_j(\sqrt{L})}(x,y)|^2\dx\Big)^{1/2}\\
&\leq
r_B^{2M-2l}V(B)w(B)^{-1/p}(2^{-(j+k)s'}r_B^{s'})\\
&~~~~~~~~~\times
\sup_{y\in B}\Big(\int_{S_k(B)}|K_{F_j(\sqrt{L})}(x,y)|^2(1+2^jd(x,y))^{2s'}\dx\Big)^{1/2}\\
&\leq
Cr_B^{2M-2l}V(B)w(B)^{-1/p}(2^{-(j+k)s'}r_B^{s'})\sup_{y\in B}\f{1}{\sqrt{V(y,2^{-j})}}\|\delta_{2^j}F_j\|_{W^q_s}\\
&\leq Cr_B^{2M-2l}V(B)w(B)^{-1/p}(2^{-(j+k)s'}r_B^{s'})\sup_{y\in
B}\f{1}{\sqrt{V(y,2^{-j})}}.
\end{aligned}
\end{equation}
For $j\geq j_0$, we have
$$\sup_{y\in B}\f{1}{V(y,2^{-j})}\leq C
\f{(2^jr_B)^n}{V(B)}.$$ This together with (\ref{eq3}) yields
\begin{equation*}
\begin{aligned}
\|F_j(\sqrt{L})b_1\|_{L^2(S_k(B))}&\leq
Cr_B^{2M-2l}V(B)w(B)^{-1/p}2^{-(j+k)s'}2^{-s'j_0}\f{(2^jr_B)^{\f{n}{2}}}{V(B)^{\f{1}{2}}}\\
&\leq
Cr_B^{2M-2l}V(B)^{1/2}w(B)^{-1/p}2^{-(j+k)s'}2^{-s'j_0}(2^jr_B)^{\f{n}{2}}
\end{aligned}
\end{equation*}
Since $w\in A_1$, we have
$$
V(B)\leq C\f{w(B)}{w(2^jB)}V(2^jB).
$$
Then,
\begin{equation*}
\begin{aligned}
\|F_j(\sqrt{L})b_1\|_{L^2(S_k(B))}&\leq
Cr_B^{2M-2l}\f{w(B)^{1/2}}{w(2^jB)^{1/2}}V(2^jB)^{1/2}w(B)^{-1/p}2^{-(j+k)s'}2^{-s'j_0}(2^jr_B)^{\f{n}{2}}\\
&\leq
Cr_B^{2M-2l}V(2^kB)^{\f{1}{2}}w(2^kB)^{-1/p}2^{-k(s'-\f{n(2-p)}{2p})}2^{(j-j_0)(\f{n}{2})-s'}.
\end{aligned}
\end{equation*}

Therefore, $$ r_B^{2l}\sum_{j\geq
j_0}\|F_j(\sqrt{L})b_1\|_{L^2(S_k(B))}\leq
C2^{-k\epsilon}r_B^{2M}V(2^kB)^{\f{1}{2}}w(2^kB)^{-1/p}
$$
with $\epsilon=s'-\f{n(2-p)}{2p}$.\\

When $j\leq j_0$, repeating the argument above with the fact that $\sup_{y\in B}\f{1}{V(y,2^{-j})}\leq C\f{1}{V(B)}$, we conclude that
$$
r_B^{2l}2^{2jM}\sum_{j< j_0}\|F_j(\sqrt{L})b_2\|_{L^2(S_k(B))}\leq
C2^{-k\epsilon}r_B^{2M}V(2^kB)^{\f{1}{2}}w(2^kB)^{-1/p}.
$$
It therefore follows that $\widetilde{a}=F(\sqrt{L})a $ is a
multiple of $(M, p, w,\epsilon)$-molecule. The proof is complete.\\

Theorem \ref{thm1} gives the boundedness on weighted Hardy spaces $\HL$
for the spectral multipliers of a non-negative self adjoint operator satisfying Gaussian upper bounds (H3).
However, there are certain operators satisfying (H1) and (H2) but not (H3) such as
general Laplace-Beltrami operators on complete Riemannian manifolds.
The following result shows that for  suitable functions $F$, the spectral operators $F(L)$ are bounded on $\HL$.

\begin{thm}\label{specmultiplierpro}
Assume that $L$ satisfies conditions $(H1)$ and $(H2)$ and $w\in A_1\cap RH_{\f{2}{2-p}}$.
Let $F$ be a bounded function defined on $(0,\infty)$ such that for
some positive number $s>\f{n(2-p)}{2p}+\f{1}{2}$ and any non-zero function $\eta\in C_c^\infty(\f{1}{2},2)$ there exists a constant $C_\eta$ such that
\begin{equation}\label{specequa}
\sup_{t>0}\|\eta(\cdot)F(t\cdot)\|_{W^2_s(\mathbb{R})}\leq C_\eta.
\end{equation}
Then the multiplier operator $F(L)$ is bounded on $\HL$.
\end{thm}
To prove Theorem \ref{specmultiplierpro} we need the following estimate in \cite{DP}.
\begin{lem}\label{spmultlem}
Let $\gamma>1/2$ and $\beta>0$. Then there exists a constant $C>0$ such that for every function $F\in W^2_{\gamma+\beta/2}$ and every function $g\in L^2(X)$ supported in the ball $B$, we have
$$
\int_{d(x,x_B)>2r_B}|F(2^j\sqrt{L})g(x)|^2\Big(\f{d(x,x_B)}{r_B}\Big)^{\beta}d\mu(x)\leq C(r_B2^j)^{-\beta}\|F\|_{W^2_{\gamma+\beta/2}}^2\|g\|_{L^2}^2
$$
for $j\in \mathbb{Z}$.
\end{lem}
\emph{Proof of Theorem \ref{specmultiplierpro}:} The proof of Theorem \ref{specmultiplierpro} can be proceeded in the same line as Theorem \ref{thm1}
by replacing Lemma \ref{lem1} by Lemma \ref{spmultlem}. We omit the details here.\\

{\bf Acknowledgement.} The authors would like to thank Lixin Yan for
helpful comments and discussions.

\newpage

\noindent The Anh Bui

\medskip

\noindent Department of Mathematics, Macquarie University, NSW 2109, Australia and
Department of Mathematics, University of Pedagogy, Ho chi Minh city, Vietnam

\smallskip

\noindent{\it E-mail:} \texttt{the.bui@mq.edu.au}, \texttt{bt\_anh80@yahoo.com}

\bigskip

\noindent Xuan Thinh Duong

\medskip

\noindent Department of Mathematics, Macquarie University, NSW 2109, Australia

\smallskip

\noindent{\it E-mail:} \texttt{xuan.duong@mq.edu.au}

\end{document}